\documentclass[12pt]{article}
\usepackage{amssymb}
\usepackage{graphicx}
\usepackage{latexsym}
\def\part#1{\frac{\partial\phantom{q}}{\partial#1}}

\newenvironment{rmk}{\begin{trivlist}\item[]{\bf Remark:} }
{\end{trivlist}}
\newenvironment{rmks}{\begin{trivlist}\item[]{\bf Remarks:} }
{\end{trivlist}}

\newenvironment{exs}{\begin{trivlist}\item[]{\bf Examples:} }
{\end{trivlist}}
\newenvironment{prf}{\begin{trivlist}\item[]{\bf Proof:} }
{\hfill $\Box$ \end{trivlist}}

\newtheorem{thm}{Theorem}

\newtheorem{prp}[thm]{Proposition}

\newcommand{\lie}[1]{\mathfrak{#1}}
\def\End{\mathop{\rm End}\nolimits}

\def\Hom{\mathop{\rm Hom}\nolimits}

\def\Sym{\mathop{\rm Sym}\nolimits}

\newcommand{\C}{\mathbf{C}}

\newcommand{\Z}{\mathbf{Z}}

\newcommand{\PP}{{\mathbf {\rm P}}}

\newcommand{\HH}{\mathbb H}
\begin{document}
\title{Stable bundles and polyvector fields}
 \author{Nigel Hitchin\\[5pt]}
\maketitle
\section{Introduction}
The moduli space of stable bundles on an algebraic curve $C$ is a much-studied object, but there are still new things to learn about it. This paper introduces one more aspect to study, and poses some conjectures about it.

Recall from \cite{NJH} that, if $M$ is such a moduli space, then its cotangent bundle $T^*M$ defines  a completely integrable  Hamiltonian system. By this we mean that there are $\dim M$ functionally independent holomorphic functions on $T^*M$ which Poisson-commute and whose common level set is an open set in an abelian variety, on which the Hamiltonian vector fields are linear. These functions are polynomial in the fibre directions and can be understood on the moduli space $M$ itself as holomorphic sections of symmetric powers $\Sym^kT$ of the tangent bundle for various values of $k$. The fact that they Poisson commute is equivalent to the statement that the symmetric tensors commute using the Schouten--Nijenhuis bracket, a natural extension of the Lie bracket on vector fields.

We introduce here a skew-symmetric version of this, identifying holomorphic sections of $\Lambda^kT$ for various values of $k$ (so-called {\it polyvector fields}) which also Schouten-commute. More precisely we note that at a smooth point of the moduli space of stable holomorphic structures on a principal $G$-bundle, where $G$ is a complex simple Lie group, the cotangent space is isomorphic to $H^0(C,{\lie g}\otimes K)$ where ${\lie g}$ denotes the adjoint bundle of Lie algebras. Given a bi-invariant differential form $\rho$ on $G$ of degree $k$, then for  $\Phi_i\in H^0(C,{\lie g}\otimes K)$, $\rho(\Phi_1,\dots,\Phi_k)$ defines a skew form with values in the line bundle $K^k$. Dually, it defines a homomorphism
$$H^1(C,K^{1-k})\rightarrow H^0(M, \Lambda^kT).$$

By analogy with the symmetric case there are three obvious questions to ask:
\begin{itemize}
\item
Is this map injective?
\item
Do these polyvector fields Schouten-commute?
\item
Is the  algebra of all polyvector fields on $M$ generated by these?
\end{itemize}
In this paper we restrict ourselves mainly to the rank one case where the only  invariant form is  $B([X,Y],Z)$ where $B$ is the Killing form, but many of our results hold in more generality.  We answer in the positive the first  question (for genus $g>2$), and show that  for general reasons the answer is yes to the second. As to the third question  the Verlinde formula shows that the answer is no, though in the final section we discuss some related issues.

 What we do show, however, is that for $g>4$, and in the case of vector bundles with coprime degree and rank, there are no polyvector fields of degree two. The vanishing in degree one is well-known \cite{NR}, \cite{NJH}, so that the first degree where existence holds is precisely  where our construction begins.  Our proof of the vanishing of $H^0(M,\Lambda^2T)$ requires another  feature of the moduli space, which was the original motivation for this research.  In \cite{NR}, the authors used a holomorphic differential on $C$ to define on the moduli space a nontrivial extension 
$$0\rightarrow T^*\rightarrow \bigoplus_{i=1}^{2g-2}{\lie g}_{x_i}\rightarrow T\rightarrow 0$$
where ${\lie g}_{x_i}$ is the restriction of the universal adjoint bundle to $M\times \{x_i\}$ and $x_i\in C$ is a zero of the differential. A considerable part of the paper consists of studying the vector bundles $E$ defined by these extensions in more detail. 

The most important point is that the extension class lies in the skew-symmetric part $H^1(M,\Lambda^2T^*)\subset H^1(M,\Hom(T,T^*))$ and realizes the known isomorphism between the space of differentials $H^0(C,K)$ and $H^1(M,\Lambda^2T^*)$. This provides an orthogonal structure on $E$ such that the subbundle $T^*$ is maximal isotropic.  We also show that these bundles have a natural Courant algebroid structure  arising from an infinite-dimensional quotient construction.

One result we need is an isomorphism $H^1(M,T\otimes T^*)\cong H^1(C,{\mathcal O})$, proved in \cite{B}. We shall see  this isomorphism being  realized as a deformation of the tangent bundle of $M$ by replacing $H^1(C,\lie{g})$ by $H^1(C,\lie{g}\otimes L)$ for a degree zero line bundle $L$. 

\section{Polyvector fields}
\subsection{The construction} \label{poly1}
We set up the basic framework in the case of a general simple Lie group $G$. Let $C$ be a compact Riemann surface and $M$ be the moduli space of stable principal $G$-bundles on $C$. At a smooth point of $M$ the tangent space $T$ is isomorphic to $H^1(C,{\lie g})$ where $\lie g$ denotes the adjoint bundle, and its dual space  $T^*$ is, by Serre duality, $H^0(C, {\lie g}\otimes K)$. We shall call sections $\Phi$ of ${\lie g}\otimes K$ Higgs fields.

Evaluation of a Higgs field at $x\in C$ defines a homomorphism from $T^*$ to ${\lie g}_x\otimes K_x$ and so a section
$$s_x\in H^0(M,{\lie g}_x\otimes T)\otimes K_x.$$

If ${\lie g}$ now denotes the universal adjoint bundle over the product $M\times C$, then varying $x$ we get a tautological section 
$$s\in H^0(M\times C,  {\lie g}\otimes(T\boxtimes K)).$$

Consider now the ring of  bi-invariant  differential forms on $G$.  This  is an exterior algebra generated by basic forms whose  degrees are given by $k_i=2m_i+1$ where $m_i$ are the exponents of the Lie algebra. For each generator $\sigma_i$ we can evaluate on the section

$$s^{\wedge k_i}\in H^0(M\times C,  \Lambda^{k_i}({\lie g}\otimes(T\boxtimes K)).$$

to obtain  
$$s_i\in H^0(M\times C, \Lambda^{k_i}T\boxtimes  K^{k_i})\cong  H^0(M, \Lambda^{k_i}T)\otimes  H^0(C,K^{k_i})$$
or equivalently by Serre duality a homomorphism 
\begin{equation}
A: H^1(C,K^{1-k_i})\rightarrow  H^0(M, \Lambda^{k_i}T).
\label{map}
\end{equation}

\begin{exs} 

\noindent 1. The simplest invariant form for any $G$ is $\sigma(X,Y,Z)=B([X,Y],Z)$ where $B$ is the Killing form. Thus  the $(5g-5)$-dimensional space $H^1(C,K^{-2})$ maps to $H^0(M, \Lambda^{3}T)$.

\noindent 2. For each point $x\in C$, evaluation of a section of $K^{k}$ at $x$ (and a trivialization of $K_x^k$) defines a linear form on $H^0(C,K^{k_i})$ and hence an element of its dual space $H^1(C,K^{1-k_i})$, so this defines a section $\sigma_x\in H^0(M, \Lambda^{k_i}T)$. 
\end{exs}

\subsection{Injectivity}
If the map $A$ in (\ref{map})   for an invariant form $\sigma$ of  degree $k$ has a non-zero kernel then there is a class $\alpha\in H^1(C,K^{1-k})$ such that  for all $G$-bundles and  Higgs fields $\Phi_1,\dots,\Phi_k$  
$$\langle \alpha, \sigma(\Phi_1,\dots,\Phi_k)\rangle=0 $$
where $\langle\,\,,\,\,\rangle$ is the Serre duality pairing. Thus for injectivity we need to show that the sections $ \sigma(\Phi_1,\dots,\Phi_k)$ generate $H^0(C,K^k)$. Here for simplicity we restrict to the rank one case.

\begin{rmk} We should make a remark here about which moduli spaces we are concerned with. The setting for the problem is the  space of stable principal $G$-bundles modulo isomorphism, but the definition of the polyvector fields only depends on the structure of the Lie algebra so it is really the adjoint group which is relevant here. In the case of a linear group the most studied moduli space is that of vector bundles of rank $n$ and degree $d$ with fixed determinant bundle. Especially important is the case where $n$ and $d$ are coprime for the moduli space then is compact and smooth and has a universal vector bundle. But it is the (singular) quotient of this by the operation of tensoring with a line bundle of order $n$ which gives the adjoint bundle moduli space. This point will become relevant in the final section. 
\end{rmk}

\begin{prp} \label{inj} If $g>2$ the map $A$ is injective for $G=SL(2)$ or $SO(3)$. 
\end{prp}

\begin{rmk} The map is not injective for $g=2$. In fact there are  two spaces, $\PP^3$ is the moduli space of bundles for even degree and  the intersection of two quadrics in $\PP^5$ for bundles of odd degree. In both cases  these are acted on trivially by the hyperelliptic involution $\tau$ on $C$ , and in particular  the action on sections of $\Lambda^3T$ is trivial. But $H^0(C,K^3)$ has both invariant and anti-invariant elements under $\tau$. 
\end{rmk}

\begin{prf} There is just one invariant form here -- the three-form $\sigma$ given by $\sigma(X,Y,Z)=B([X,Y],Z)$. 

\noindent 1.  We begin with the even degree case,  and we may consider a class  to be represented by a rank $2$ vector bundle $E$ with $\Lambda^2 E$ trivial. Consider first a non-trivial extension of degree zero line bundles 
\begin{equation}
0\rightarrow L\rightarrow E\rightarrow L^*\rightarrow 0
\label{ext}
\end{equation}
defined by $\alpha\in H^1(C,L^2)$. Take  a point $x\in C$, and let $t$ denote the tautological section of ${\mathcal O}(x)$. The homomorphism $L^*(-x)\rightarrow L^*$ defined by the product with $t$ lifts to $E$ if the class $\alpha t\in H^1(C,L^2(x))$ vanishes. The long exact sequence of 
$$0\rightarrow {\mathcal O}_C(L^2)\stackrel {t}\rightarrow {\mathcal O}_C(L^2(x))\rightarrow   {\mathcal O}_x(L^2(x))\rightarrow 0$$
gives
$$\rightarrow H^0(C,L^2(x))\rightarrow L^2(x)\vert_x\rightarrow H^1(C,L^2)\stackrel{t}\rightarrow  H^1(C,L^2(x))\rightarrow$$
so if $H^0(C,L^2(x))=0$ there is a unique $\alpha$ with this property -- the image of a vector in $L^2(x)\vert_x$ under the connecting homomorphism. Moreover from the exact sequence of (\ref{ext}) the lift is then unique. 

The lift defines a section of $\Hom(L^*(-x),E)=EL(x)$. This is an inclusion unless it vanishes at $x$ but if that were so,  it would come from a section of $EL$ and in the long exact sequence of
$$0\rightarrow L^2\rightarrow EL\rightarrow {\mathcal O}\rightarrow 0$$
we see that the generator of $H^0(C,{\mathcal O})$ maps to $\alpha\in H^1(C,L^2)$ so if $L^2$ is non-trivial, $H^0(C,EL)=0$.

Hence if $H^0(C,L^2(x))=0$ (which implies $L^2$ is non-trivial), the lift of  $t: L^*(-x)\rightarrow L^*$ to $E$ gives another expression of $E$ as an extension
$$0\rightarrow L^*(-x)\rightarrow E\rightarrow L(x)\rightarrow 0.$$

 If $H^0(C,L^2(x))\ne 0$ then there is a point $y$ such that the divisor class $[L^2]\sim y-x$. This defines a two-dimensional subvariety of the Jacobian  and hence for $g>2$ a generic line bundle $L$ has the property that for all $x$, $H^0(C,L^2(x))=0$.

 From \cite{NR1} a generic element in the $(g+1)$-dimensional space $H^1(C,L^{-2}(-2x))$ defines an extension as above which is a  stable bundle.   Thus, as  we  vary $L$ and $x$ and the extension class, $E$ belongs to a family whose generic member is stable. Moreover, although $E$ itself is not stable it is simple, i.e. it has no non-scalar  endomorphisms. This means that the rank of $H^0(C,{\lie g}\otimes K)=3g-3$ for all bundles in the family. 

We shall show that, varying $L$ and $x$, $\sigma(\Phi_1,\Phi_2,\Phi_3)$ generates $H^0(C,K^3)$ and hence will do so in  a generic family of stable bundles.
\vskip .5cm
\noindent 2. We now need to determine the Higgs fields for $E$. The adjoint bundle ${\lie g}$ is $\End_0 E$, the bundle of trace zero endomorphisms. We have an exact sequence
$$0\rightarrow E\otimes LK\rightarrow \End_0 E\otimes K \rightarrow L^{-2}K\rightarrow 0$$
and so a section $s$ of $L^{-2}K$ lifts to a Higgs field $\Phi_1$ if it maps  in the long exact sequence to zero in $H^1(C,E\otimes LK)$. Now $H^1(C,L^2K)$ is dual to $H^0(C,L^{-2})$ which vanishes if $L^2$ is non-trivial which means from the long exact sequence of (\ref{ext}) that $H^1(C,E\otimes LK)\cong H^1(C,K)$ and so $s\in H^0(C,L^{-2}K)$ extends if its product with the extension class $\alpha\in H^1(C,L^{2})$ vanishes. Choosing the class as above, this means that $s(x)=0$. 

Now let $\Phi_1$ be a lift of $s$.  In the exact sequence 
$$0\rightarrow L^2K\rightarrow E\otimes LK\rightarrow K \rightarrow 0$$
since $H^1(C,L^2K)=0$ the map $H^0(C,E\otimes LK)\rightarrow H^0(C,K)$ is surjective, so given a section $t$ of $K$ we can find $\Phi_2$ a section of $E\otimes LK\subset  \End_0 E\otimes K$ such that $\Phi_2$ maps to $t$. Now choose $\Phi_3$ to be any section $u$ of $L^2K\subset E\otimes LK\subset  \End_0 E\otimes K$. The $(3g-3)$-dimensional space of Higgs fields can now be seen to be constructed from $s$, in the $(g-2)$-dimensional subspace of  $H^0(C,L^{-2}K)$ consisting of sections that vanish at $x$, a choice of $t$ in the $g$-dimensional space of differentials, and an arbitrary section $u$ in the $(g-1)$-dimensional space of sections of $L^{2}K$.

It is easy to see then that $\sigma(\Phi_1,\Phi_2,\Phi_3)\in H^0(C,K^3)$ is a multiple of $stu$. We shall show that we can generate all sections of $K^3$ vanishing at $x$ this way, for $L$ generic. 
\vskip .5cm
\noindent 3. We use the ``base-point free pencil trick" of \cite{ACGH}: let  $U$ be a line bundle with sections $s_1,s_2$ having no common zeros and let $F$ be a vector bundle.  Then the kernel of the map
$$\C^2\otimes H^0(C,F)\rightarrow H^0(C,F\otimes U)$$
defined by $(t_1,t_2)\mapsto s_1t_1+s_2t_2$ is isomorphic to $H^0(C,F\otimes U^*)$. Indeed, if $s_1t_1=-s_2t_2$ and $s_1,s_2$ have no common zeros then $t_1=us_2, t_2=-us_1$ for $u$ a section of $F\otimes U^*$.

The bundle $L^2K$ has a basepoint $x$ if $H^1(C,L^2K(-x))\ne 0$, or by Serre duality if  $H^0(C,L^{-2}(x))\ne 0$. As above,  considering $[L^2]\sim y-x$, if  $g>2$ it follows that  for generic $L$,  $L^2K$ has no basepoint. Consider the sequence
$$0\rightarrow {\mathcal O}_C L^2K(-x-y)\rightarrow  {\mathcal O}_C L^2K(-x)\rightarrow  {\mathcal O}_y L^2K(-x)\rightarrow 0.$$
In the long exact sequence we see that  if $H^1(C,L^2K(-x-y))\rightarrow H^1(C,L^2K(-x))$ is injective, sections of $L^2K$ vanishing at $x$ do not all vanish at a given point $y$. This injectivity condition is equivalent to the map $H^0(C,L^{-2}(x))\rightarrow H^0(C,L^{-2}(x+y))$ being surjective. But  if $L^{-2}(x+y)$ has a section then $[L^2]\sim x+y-u-v$. Thus if the genus $g>4$ then for  generic $L$ there are no sections. By Riemann-Roch $\dim H^0(C,L^2K)=g-1$; this system has no basepoint and separates points hence  the map $C\rightarrow \PP^{g-2}$ is injective and we can then use general position arguments.

Take a general divisor $x_1+\dots+x_{2g-2}$ of $L^{2}K$. By general position the space of sections vanishing on the first $(g-3)$ of these points  is a base-point free two-dimensional space: sections of $L^2K(-D)$ where $D=x_1+\dots+x_{g-3}$. Taking $F=K$ in the ``base-point free pencil trick", we have, by taking products with $s_1$ and $s_2$, a subspace of $H^0(C,L^2K^2)$ of dimension 
$2g-\dim H^0(C, L^{-2}(D))$. But by Riemann-Roch 
$$\dim H^0(C, L^{-2}(D))-\dim H^1(C, L^{-2}(D))=g-3+1-g=-2$$
and by Serre duality $\dim H^1(C, L^{-2}(D))=\dim H^0(C,L^2K(-D))=2$. Hence we obtain $H^0(C, L^{-2}(D))=0$.
This means that $s_2,s_2$ generate with sections of $K$ a $2g$-dimensional subspace of $H^0(C,L^2K^2)$, which by Riemann-Roch has dimension $(3g-3)$.

However, by general position, for each point $x_i $ there is a section $s_i$ of $L^2K$ which vanishes at all points $x_1,\dots, x_{g-3}$  except $x_i$.  Multiplying these by sections of $K$ gives a complementary $(g-3)$-dimensional subspace and hence in total $2g+(g-3)=3g-3$ linearly independent sections. Hence sections of $L^2K^2$ are generated by sections of $L^2K$ and $K$.

Given  $x \in C$ and  $2g-3$ general points $x_1,x_2,\dots,x_{2g-3}$ there is a section $q$ of $K^3$ vanishing at these points since $K^3$ defines an embedding  $C\subset \PP^{5g-6}$. But since $2g-3\ge g$, $2g-3$  general points form the divisor for a generic line bundle $L^{-2}K(-x)$. Hence the divisor of $q$ is of the form  $x+D_1+D_2$ where $D_1$ is the divisor of a section of $L^{-2}K(-x)$ and $D_2$ of $L^2K^2$. Using the above result about sections of $L^2K^2$ we see that  sections of $L^{-2}K(-x),K,L^2K$ generate $H^0(C,K^3(-x))$. 

 Varying $x$, by genericity we can generate all sections of $K^3$  from extensions $E$ and hence also from stable bundles.
\vskip .5cm
There remain the cases $g=3,4$. For $g=4$ a generic degree $2g-2=6$ line bundle maps $C$ birationally to a singular sextic curve in $\PP^2$ and the genericity theorem holds here. For genus $g=3$, the image of  $H^0(C,L^2K)\otimes H^0(C,K)$ has dimension $2\times 3=6=3g-3$ by the base-point trick. 
\vskip .5cm
The case of odd degree can be considered as the study of the moduli space of rank $2$ vector bundles $E$ with $\Lambda^2E\cong {\mathcal O}(-y)$.  Here, from \cite{NR1} each non-trivial extension 
\begin{equation}
0\rightarrow L\rightarrow E\rightarrow L^*(-y)\rightarrow 0
\label{extodd}
\end{equation}
is stable (and indeed each stable bundle arises this way). The extension is defined by a class $\alpha\in H^1(C,L^2(y))$ as before and $s\in H^0(C,L^{-2}K(-y))$ lifts to a Higgs field $\Phi_1$ if its product with $\alpha$ vanishes. For some point $x$ we can take as above an extension where, for $g>4$,  this condition is that $s(x)=0$. We now need to prove that any section of $K^3$ vanishing at $x$ is generated by sections of the $g$-dimensional spaces $H^0(C,L^2K(y)), H^0(C,K)$ and the $(g-3)$-dimensional space $H^0(C,L^{-2}K(-x-y))$ for some $y$. The argument proceeds as before when $L$ is chosen so that sections of $L^{-2}K$ map $C$ birationally to its image. For genus $3$, $L^2K(y)$ defines an embedding for generic $L$, so sections of $L^2K(y)$ and $K$ generate sections of $L^2K^2(y)$.
\end{prf}

\begin{rmk} To consider the case of general $G$ we could at this point use our rank one starting point and take the homomorphism from $SL(2)$ to $G$ given by the principal three-dimensional subgroup \cite{K} to define a $G$-bundle. This breaks up the Lie algebra ${\lie g}$ into irreducible representations of $SO(3)$  whose dimensions are precisely the degrees of the generators of the algebra of invariant differential forms on $G$. It seems reasonable to conjecture that the restriction of a generating form of degree $(2m+1)$ to the corresponding subspace of the same dimension is non-zero, but this seems not to have been proved, except for $m=1$ where it is clear. If it were true then a direct generalization of the above  would give injectivity in general, though  there may well be other means of achieving this.  As it stands we can assert that $H^1(C,K^{-2})\rightarrow H^0(M,\Lambda^3T)$ is injective for all $G$.
\end{rmk}

\subsection{The Schouten-Nijenhuis bracket} \label{Sch}

If $A,B$ are sections of $\Lambda^p T$ and $\Lambda^q T$ respectively on any smooth manifold then one can form the {\it Schouten-Nijenhuis bracket} $[A,B]$ which is  a section of  $\Lambda^{p+q-1}T$, generalizing the Lie bracket of two vector fields. It has the basic properties:

\begin{itemize}
\item
For each vector field $X$, $[X,A]={\mathcal L}_XA$
\item
$[A,B]=-(-1)^{(p-1)(q-1)}[B,A]$
\item
$[A,B\wedge C]=[A,B]\wedge C+(-1)^{(p-1)q}B\wedge [A,C]$
\end{itemize}

\begin{rmk} There is a similar bracket on sections of $\Sym^p T$ and $\Sym^q T$ which corresponds to the Poisson bracket of the corresponding functions on the total space of the cotangent bundle $T^*$ with respect to the  canonical symplectic form.
\end{rmk}

On a complex manifold, the Schouten-Nijenhuis bracket on the sheaf of holomorphic polyvector fields extends to give a Gerstenhaber algebra structure on $H^*(M,\Lambda^*T)$. We shall show here that the global polyvector fields just constructed commute with respect to this bracket. 

 We adopt an infinite-dimensional viewpoint which can be made rigorous in a standard way by using Banach manifolds and slice theorems. Consider the moduli space $M$ of stable bundles as the quotient of an open set in the  space ${\mathcal A}$ of all $\bar\partial$-operators $\bar\partial_A$ on a fixed $C^{\infty}$ bundle by the group ${\mathcal G}$ of complex gauge transformations. The space ${\mathcal A}$ is an infinite-dimensional affine space with translation group $\Omega^{0,1}(C,{\lie g})$. The cotangent space at any point is formally $\Omega^0(C,{\lie g}\otimes K)=\Omega^{1,0}(C,{\lie g})$ using the pairing for 
 $a\in \Omega^{0,1}(C,{\lie g})$ and $\Phi \in \Omega^0(C,{\lie g}\otimes K)$,
 $$\int_CB(\Phi,a)$$
 where $B$ is the Killing form.
 
 Take $\alpha\in \Omega^{01}(C,K^{1-k})$ representing a class in $H^1(C,K^{1-k})$ and an invariant $k$-form $\sigma$ on ${\lie g}$. Then define a  $k$-vector field $S$ on ${\mathcal A}$ by evaluating on cotangent vectors $\Phi_i\in \Omega^0(C,{\lie g}\otimes K)$:
 $$S(\Phi_1,\dots,\Phi_k)=\int_C\sigma(\Phi_1,\dots,\Phi_k) \alpha.$$
 Since $\alpha$ is independent  of the operator $\bar\partial_A$  such a polyvector field on ${\mathcal A}$ is translation invariant (has ``constant coefficients"), so any two Schouten-commute. 

But $S$ is gauge-invariant because $\sigma$ is invariant, so under the derivative of the quotient map from the open set of stable points in ${\mathcal A}$ to $M$, 
$$\Lambda^kT_A{\mathcal A}\rightarrow \Lambda^kT_{[A]}M$$
the image is independent of the representative point $A$, and so defines a polyvector field $\bar S$ on $M$. Note that an invariant polyvector field $S$ is not the same as a polyvector field $\bar S$ on the quotient but $\bar S$ is defined by evaluating on 1-forms which are pulled back. In our case these are {\it holomorphic} sections $\Phi_i$, and then $\sigma(\Phi_1,\dots,\Phi_k)$ is a holomorphic section of $K^k$, and by Stokes' theorem only the Dolbeault cohomology class of $\alpha$ contributes  in the definition.

We shall show in Section \ref{vanish} that there are no holomorphic bivector fields, sections of $\Lambda^2T$, for $g>4$ on the moduli space of stable vector bundles when the rank and degree are coprime. Since we have seen above that there always exist  non-trivial holomorphic trivector fields, this is some information towards answering the third question in the Introduction. The result most probably extends to other groups but we shall use theorems in the literature which relate to this particularly familiar case. 

Our approach revisits a vector bundle on the moduli space first introduced by Narasimhan and Ramanan \cite{NR}, but where we observe some extra features.

\section{Orthogonal bundles on the moduli space}

\subsection{Courant algebroids}\label{cou}

We need the notion of a holomorphic (exact) Courant algebroid. This is  a vector bundle $E$ given as an extension 
$$0\rightarrow T^*\rightarrow E\stackrel{\pi} \rightarrow T\rightarrow 0$$
with the following properties. 
\begin{itemize}
\item
$E$ has an orthogonal structure -- a nondegenerate symmetric form $(\,\,,\,\,)$ such that $T^*\subset E$ is isotropic.
\item
For local sections $u,v$ there is another local section $[u,v]$, skew-symmetric in $u,v$, such that  :

\noindent (i) if $f$ is a local function $[u,fv]=f[u,v]+(\pi(u)f)v-(u,v)df$

\noindent (ii) $\pi(u)(v,w)=([u,v]+d(u,v),w)+(v,[u,w]+d(u,w))$

\noindent (iii) $[u,[v,w]]+[v,[w,u]]+[w,[u,v]]=d(([u,v],w)+([w,u],v)+([v,w],u))/3$
\end{itemize}

The standard example is $T\oplus T^*$ with symmetric form
$$(X+\xi,X+\xi)=i_X\xi$$ and  bracket
$$[X+\xi, Y+\eta]=[X,Y]+\mathcal{L}_{X}\eta -\mathcal{L}_{Y}\xi -\frac{1}{2}
  d(i_{X}\eta  -i_{Y}\xi).$$
  
  There is a natural quotient construction for Courant algebroids (see \cite{Quot}) which we describe next. Suppose $M$ is a manifold with a free proper action of a Lie group $G$, and suppose that there is a lifted action on $E$ preserving all the structure, in particular being compatible with $\pi: E\rightarrow T$ and the natural action on $T$.    The derivative of the group action defines a Lie algebra homomorphism $a\mapsto X_a$ from ${\lie g}$ to vector fields, sections of $T$, and we ask for an equivariant {\it extended action} which is a linear map $e$ from ${\lie g}$ to sections of $E$ such that:
  \begin{itemize}
  \item
   for $a,b\in {\lie g}$ we have $[e(a),e(b)]=e([a,b])$
   \item
     $(e(a),e(a))=0$ 
 \item
  $\pi(e(a))=X_a$
  \end{itemize}
  
  Given this data, $e({\lie g})$ generates a trivial subbundle  $F\subset E$ of rank $\dim G$. It is  isotropic by the second condition, so $F\subset F^{\perp}$, and $F^{\perp}/F$ inherits a nondegenerate symmetric form.  The latter is a $G$-invariant bundle of rank $(2\dim M-2\dim G)=2\dim (M/G)$. By $G$-invariance it descends to a bundle $\bar E$ on $M/G$. The $G$-invariant sections of $F^{\perp}/F$ are by definition the sections of $\bar E$ on $M/G$ and the bracket on $G$-invariant sections defines a bracket on sections of $\bar E$. 
  
  Now $\pi(F)$ is the tangent bundle along the fibres of $M\rightarrow M/G$, so $\pi$ induces a map from $F^{\perp}/F$ to $T(M/G)$ and one can easily deduce that  $\bar E$ is  a Courant algebroid over $M/G$. 
  
  \subsection{A family of Courant algebroids}
  We shall give here an infinite-dimensional example of the above construction to produce (for general $G$) a family of Courant algebroids over the moduli space of stable $G$-bundles.
  
  As in Section \ref{Sch} let ${\mathcal A}$ denote the infinite-dimensional  space of all holomorphic structures on a fixed principal $G$-bundle.  It is acted on by the group ${\mathcal G}$ of complex gauge transformations, and the quotient of the open set of stable holomorphic structures by ${\mathcal G}$ is the finite-dimensional moduli space of dimension $\dim G(g-1)$. So ${\mathcal A}$ is our manifold with ${\mathcal G}$-action and we are going to define an extended action on the trivial Courant algebroid $T\oplus T^*$. For this we consider  as above the cotangent space to be $\Omega^{1,0}(C,{\lie g})$.
    
 To define an extended action we choose  a holomorphic $1$-form  $\theta\in H^0(C,K)$ and define, for $\psi$ in the Lie algebra $\Omega^0(C,{\lie g})$ of ${\mathcal G}$,
 $$e(\psi)(a)=(\bar\partial_A \psi, \psi\theta)\in  \Omega^{0,1}(C,{\lie g})\oplus \Omega^{1,0}(C,{\lie g}).$$
 
 We check the isotropy condition:
 $$(e(\psi),e(\psi))=\int_C B( \psi\theta,\bar\partial_A \psi)=\frac{1}{2}\int_C\bar\partial (\theta B(\psi,\psi))=0$$
 since $\theta$ is holomorphic.  
 
 To check the bracket condition $[e(\psi),e(\psi')]=e([\psi,\psi'])$ note that $\psi\theta$ is independent of $A$ and is thus a translation-invariant 1-form on ${\mathcal A}$ and hence is closed. Thus using ${\mathcal L}_X=di_X+i_Xd$ we have, for $\xi=\psi\theta,\eta=\psi'\theta$, $X=\bar\partial_A\psi,Y=\bar\partial_A\psi'$
$$ \mathcal{L}_{X}\eta -\mathcal{L}_{Y}\xi -\frac{1}{2}
  d(i_{X}\eta  -i_{Y}\xi)=\frac{1}{2}
  d(i_{X}\eta  -i_{Y}\xi).$$
  Now
  $$i_{X}\eta  -i_{Y}\xi=\int_C B(\psi'\theta,\bar\partial_A\psi)-B(\psi\theta,\bar\partial_A\psi')$$
 and $d(i_{X}\eta  -i_{Y}\xi)$ evaluated on $a\in \Omega^{0,1}(C,{\lie g})$ is
 $$\int_C B(\psi'\theta,[a,\psi])-B(\psi\theta,[a,\psi'])=-\int_C \theta B([\psi',\psi],a)-\theta B([\psi,\psi'],a)=2\int_C\theta B([\psi,\psi'],a).$$
 
 This does then define an extended action and we can produce a quotient Courant algebroid as in Section \ref{cou}. In our case the space $F$ generated by the Lie algebra of ${\mathcal G}$ consists of the subspace 
 $$B^1=\{(\bar\partial_A \psi, \psi\theta)\in  \Omega^{0,1}(C,{\lie g})\oplus \Omega^{1,0}(C,{\lie g})\}$$ and $F^{\perp}$ is the space of pairs $(a,\Phi)\in  \Omega^{0,1}(C,{\lie g})\oplus \Omega^{1,0}(C,{\lie g})$ such that
 $$\int_C B(\Phi,\bar\partial_A \psi)+B(\psi\theta,a)=0$$
 for all $\psi\in \Omega^0(C,{\lie g})$. By integration by parts this is 
 $$Z^1=\{(a, \Phi)\in  \Omega^{0,1}(C,{\lie g})\oplus \Omega^{1,0}(C,{\lie g}): \bar\partial_A\Phi=a\theta\}.$$
Hence $F^{\perp}/F=Z^1/B^1$ is the first cohomology group  of the complex
$$ \Omega^{0}(C,{\lie g})\stackrel{\bar\partial +\theta}\rightarrow \Omega^{0,1}(C,{\lie g})\oplus \Omega^{1,0}(C,{\lie g})\stackrel{\bar\partial +\theta}\rightarrow  \Omega^{1,1}(C,{\lie g})$$
or equivalently the hypercohomology  $\HH^1(C,{\lie g})$ of the  short complex  of sheaves
 $$ {\mathcal O}({\lie g})\stackrel{\theta}\rightarrow {\mathcal O}({\lie g}\otimes K).$$
 From the first hypercohomology spectral sequence we have an exact sequence
  $$H^0(C,{\lie g})\rightarrow H^0(C,{\lie g}\otimes K)\rightarrow \HH^1(C,{\lie g})\rightarrow H^1(C,{\lie g})\rightarrow  H^1(C,{\lie g}\otimes K)$$
 which for stable bundles gives us the expected extension
 $$0\rightarrow T^*\rightarrow \HH^1(C,{\lie g})\rightarrow T\rightarrow  0.$$
 For the second sequence, if ${\mathcal Q}$ is the quotient sheaf
  $$0\rightarrow  {\mathcal O}({\lie g})\stackrel{\theta}\rightarrow {\mathcal O}({\lie g}\otimes K)\rightarrow {\mathcal Q}\rightarrow 0$$
  we have
 $$0\rightarrow \HH^1(C,{\lie g})\stackrel{\cong}\rightarrow H^0(C,{\mathcal Q})\rightarrow 0.$$
But ${\mathcal Q}$ is supported on the zero-set of the differential $\theta$. So for generic $\theta$ with simple zeros  $x_1,\dots,x_{2g-2}$ we have an isomorphism  from $\HH^1(C,{\lie g})$ to
$$\bigoplus_{i=1}^{2g-2}({\lie g}\otimes K)_{x_i}.$$
Denoting  by ${\lie g}_x$ the universal adjoint bundle restricted to $M\times \{x\}$ we find that the Courant algebroid $E$ on $M$ produced by our quotient construction is a direct sum of bundles 
\begin{equation}
E\cong \bigoplus_{i=1}^{2g-2}{\lie g}_{x_i}\otimes K_{x_i}.
\label{sum}
\end{equation}
\begin{rmks}

\noindent 1. This vector bundle and its description as an extension appeared in the paper \cite{NR}. It is the simplest way to see that the total Pontryagin class of $M$ (the total Chern class of $T\oplus T^*$) is of the form $p(T)=c({\lie g})^{2g-2}$. Neither the symmetric form nor the Courant bracket played a role in its initial introduction.

\noindent 2. The extension $0\rightarrow T^*\rightarrow E\rightarrow T\rightarrow 0$ defines a class in $H^1(M,T^*\otimes T^*)$ but the orthogonal structure, and the fact that $T^*$ is isotropic, tells us that the class lies in   $H^1(M,\Lambda^2T^*).$ Each such extension depended on a choice of differential $\theta$ so we have a natural homomorphism $H^0(C,K)\rightarrow H^1(M,\Lambda^2 T^*)$. For vector bundles this is an isomorphism -- in fact from \cite{NR} deformations of the point $x\in C$ give non-trivial deformations of ${\lie g}_{x}$ and hence from (\ref{sum}) non-trivial deformations of $E$, in particular non-trivial extension classes, so the map is injective; but  both spaces are $g$-dimensional. The more usual description of this isomorphism is the dual one -- the intermediate Jacobian of $M$ is isomorphic to the Jacobian of $C$ \cite{NR} hence $H^1(C,{\mathcal O})\cong H^2(M,\Lambda^1T^*)$.
\end{rmks}

\subsection{The orthogonal structure}
Let $A$ be a holomorphic structure on the principal bundle and $e\in E_{[A]}$ a vector in the fibre of $E$ over $[A]\in M$. Then $e$ is represented by $(a,\Phi)\in \Omega^{0,1}(C,{\lie g})\oplus \Omega^{1,0}(C,{\lie g})$ where
$$\bar\partial_A\Phi=a\theta.$$
The inner product is defined by 
$$(e,e)=\int_C B(\Phi,a).$$
Surround each zero of $\theta$ by a small disc and let $C'$ be the complement of these discs, then $a=\theta^{-1}\bar\partial_A\phi$ is smooth on $C'$ and 
$$\int_{C'}  B(\Phi,a)=\int_{C'} \frac{1}{\theta} B(\Phi,\bar\partial_A\Phi)=\frac{1}{2}\int_{C'}\bar\partial( \frac{1}{\theta} B(\Phi,\Phi))=0.$$
It follows directly on shrinking the discs that, for simple zeros of $\theta$, the orthogonal structure  is
\begin{equation}
(e,e)=\pi i \sum_{i=1}^{2g-2} \frac{B(\Phi,\Phi)_{x_i}}{\theta'(x_i)}
\label{ee}
\end{equation}
where $\theta'(x_i)\in K^2_{x_i}$ is the derivative of $\theta$ at its zero $x_i$. 

\begin{rmks}

\noindent 1. Note from this description of the inner product that  the decomposition of $E$ in (\ref{sum}) is an orthogonal one. 

\noindent 2. Note also that if  $\Phi$ is holomorphic then $B(\Phi,\Phi)/\theta$ is a meromorphic differential and the sum of its residues is therefore zero. Hence from (\ref{ee}) $T^*\subset E$ is maximally isotropic.
\end{rmks}

We can generalize the above  by replacing $\theta$ by  a section $s$ of $KL^2$ where $L$ is a line bundle of degree zero and considering the hypercohomology of 
 $$ {\mathcal O}({\lie g}\otimes L^*)\stackrel{s}\rightarrow {\mathcal O}({\lie g}\otimes KL).$$
The quadratic form is defined in the same way as (\ref{ee}), and $H^0(C,{\lie g}\otimes KL)$ is still isotropic but we have lost the Courant bracket.  

What we obtain this way is a hypercohomology group 
$$0\rightarrow T_L^*\rightarrow \HH^1(C,{\lie g\otimes L^*})\rightarrow T_L\rightarrow  0.$$
where $T_L= H^1(C,{\lie{g}}\otimes L^*)$. In particular, varying over the moduli space, we see that each line bundle $L$ of degree zero defines a deformation $T_L$ of the tangent bundle. 

To summarize, for each effective divisor $D$ of degree $2g-2$ we have produced an orthogonal bundle $E_D$ with the following properties 
\begin{itemize}
\item
$E_D$ has an orthogonal structure
\item
there is an exact sequence of vector bundles $0\rightarrow T_L^*\rightarrow E_D\rightarrow T_L\rightarrow 0$ where $KL^2$ is the line bundle defined by $D$
\item
$T_L^*$ is a maximal isotropic subbundle
\item
when $L$ is trivial, $E_D$ has the structure of a holomorphic Courant algebroid
\end{itemize}

\section{A vanishing theorem}\label{vanish}

We shall use the bundles $E_D$ to prove the following vanishing theorem:

\begin{thm} Let $M$ be the moduli space of rank $n$, degree $d$ bundles of  fixed determinant, with $n,d$ coprime,  over a curve of genus $g>4$. Then $H^0(M,\Lambda^2T)=0$.
\end{thm}
\begin{prf}
 We return to the situation of a 1-form $\theta$ defining an extension
$$0\rightarrow T^*\rightarrow E\rightarrow T\rightarrow 0$$
There is  an induced sequence of vector bundles 
\begin{equation}
0\rightarrow A\rightarrow \Lambda^2 E\rightarrow \Lambda^2T\rightarrow 0
\label{wedge}
\end{equation}
which we shall use to approach $\Lambda^2T$. 
Here $A$ is the bundle of Lie subalgebras preserving $T^*$ and is itself  an extension 
\begin{equation}
0\rightarrow \Lambda^2T^*\rightarrow A\rightarrow T\otimes T^*\rightarrow 0.
\label{Aext}
\end{equation}

Consider the bundle $\Lambda^2E$. From (\ref{sum}) we have 
\begin{equation}
\Lambda^2E\cong \bigoplus_{i<j}({\lie g}_{x_i}\otimes {\lie g}_{x_j})\oplus \bigoplus_{i}\Lambda^2 {\lie g}_{x_i}
\label{L2E}
\end{equation}
The coprime condition means that there is a universal vector bundle. In \cite{LN} the authors show that vector bundles $U_x$ on $M$ coming from this universal bundle are stable and isomorphic if and only if $x=y$. If ${\lie g}_{x}\otimes {\lie g}_{y}=\End_0 U_x\otimes \End_0 U_y$ has a holomorphic section then  by stability this is covariant constant with respect to the connection defined by the  Hermitian-Einstein connections on $U_x$ and $U_y$. This connection has holonomy $U(n)\cdot U(n)$ which means in particular that  the section defines an algebra homomorphism from $\End_0 U_x$ to $\End_0 U_y$. By stability this is an isomorphism which means that $U_x\cong L\otimes U_y$ for a line bundle $L$. 
But the Picard variety of $M$ is $\Z$ and $c_1(U_x)=c_1(U_y)$ so $U_x\cong U_y$ and $x=y$. It follows that $H^0(M,{\lie g}_{x_i}\otimes {\lie g}_{x_j})=0$ if $i\ne j$. If $x_i=x_j$ then we similarly deduce that the only holomorphic section of  ${\lie g}_{x_i}\otimes {\lie g}_{x_i}$ is defined by the Killing form $B$, which is symmetric and hence $H^0(M,\Lambda^2{\lie g}_{x_i})=0$. From (\ref{L2E}) we see that 
$$H^0(M,\Lambda^2 E)=0.$$

Since $E$ has an orthogonal structure, $\Lambda^2E$ is isomorphic to the bundle of skew-adjoint transformations of $E$ and the derivative of any family of deformations of $E$ as an orthogonal bundle defines an element of  $H^1(C,\Lambda^2E)$. 

But we saw in the previous section that any holomorphic section with divisor $D$ of a line bundle $KL^2$ (i.e. any bundle  of degree $(2g-2)$)  defines an extension 
$$0\rightarrow T_L^*\rightarrow E_D\rightarrow T_L\rightarrow 0$$
with an orthogonal structure  such that $T_L^*$ is isotropic. We therefore have a family of extensions  defined by a $2^{2g}$-fold covering (the choice of the line bundle $L$) of the symmetric product $S^{2g-2}C$ all of which have orthogonal structures.  So we have an effectively parametrized $(2g-2)$-dimensional family of bundles deforming $E$.  Each of these bundles has an orthogonal structure  so the tangent space to the family is a distinguished $(2g-2)$-dimensional subspace of $H^1(M,\Lambda^2 E)$. But this family also comes with a distinguished  maximal isotropic subbundle, so  this subspace is the image of  a $(2g-2)$-dimensional subspace $V\subseteq H^1(M,A)$.

Now consider the long exact sequence for (\ref{Aext}) 
$$\rightarrow H^0(M,T\otimes T^ *)\stackrel{h}\rightarrow H^1(M,\Lambda^2T^*)\rightarrow H^1(M,A)\stackrel{p}\rightarrow H^1(M,T\otimes T^ *)\rightarrow H^2(M,\Lambda^2T^*)\rightarrow$$

 From  \cite{B} for $g>4$ $H^0(M,T\otimes T^*)$ consists of  multiples of the identity. The  homomorphism $h$ is just the extension class defining $E$ in $H^1(M,\Lambda^2T^*)$ applied to the identity and so is injective. We know that $H^1(M,\Lambda^2T^*)\cong H^0(C,K)$, hence from the exact sequence the kernel of $p$ has dimension $(g-1)$.

 Now a deformation of $E$, as a bundle with distinguished subbundle, defines a deformation of the subbundle.   The map $p$ in the exact sequence is its derivative.  Our $(2g-2)$-dimensional family of deformations of $E$ is parametrized by an effective degree $(2g-2)$ divisor $D$ and defines the deformation $T_L$ of the tangent bundle, where the divisor class of $D$ is $K+2L$. This  map factors through the Abel-Jacobi map  $u: S^{2g-2}C\rightarrow J(C)$ at the divisor of $\theta$, and so $p$, restricted to the subspace $V\subseteq H^1(M,A)$, factors through the derivative of $u$.
 
  Writing the map $u$ as 
 $$u_{\alpha}=\sum_{i=1}^{2g-2}\int_{x_0}^{x_i}\omega_{\alpha}$$
 for a basis $\{\omega_{\alpha}\}$ of differentials we see that the image of its derivative is the $(g-1)$-dimensional subspace of $H^1(C,{\mathcal O})$ annihilated by $\theta\in H^0(C,K)=H^1(C,{\mathcal O})^*$. The kernel of $p$ restricted to $V$ is thus $(g-1)$-dimensional and hence coincides with the full kernel of $p$. Hence $p(V)\subset H^1(M,T\otimes T^*)$ is $(g-1)$-dimensional.  From  \cite{B}  $H^1(M,T\otimes T^*)\cong H^1(C,{\mathcal O})$ and thus has dimension $g$. We deduce that either $p:H^1(M,A)\rightarrow H^1(M,T\otimes T^ *)$ is surjective, and then $\dim H^1(M,A)=(2g-2)+1$ or $p$ maps to a $(g-1)$-dimensional space which means that $V=H^1(M,A)$ and  $\dim H^1(M,A)=2g-2.$
 \vskip .25cm
 
 Now consider the long exact sequence for (\ref{wedge}) 
$$\rightarrow H^0(M,\Lambda^2 E)\rightarrow H^0(M, \Lambda^2T)\rightarrow H^1(M,A)\rightarrow H^1(M,\Lambda^2E)\rightarrow.$$
If $V=H^1(M,A)$ then knowing that $V$ maps injectively to $H^1(M,\Lambda^2E)$ and $H^0(M,\Lambda^2E)=0$, we have the required result $H^0(M,\Lambda^2 T)=0$. The other alternative is that $\dim  H^1(M,A)=2g-1$ in which case $\dim H^0(M,\Lambda^2 T)\le 1$
\vskip .25cm
We now use the  exact sequence obtained by tensoring $E$ with $T$
 $$0\rightarrow T\otimes T^*\rightarrow T\otimes E\rightarrow T\otimes T\rightarrow 0$$
 to yield the exact cohomology sequence
 $$0\rightarrow H^0(M,T\otimes T^*)\rightarrow H^0(M,T\otimes E)\rightarrow H^0(M,T\otimes T)\rightarrow H^1(M,T\otimes T^*)\rightarrow \cdots$$
Now 
$$E\cong \bigoplus_{i=1}^{2g-2}{\lie g}_{x_i}$$
so for the term $H^0(M,T\otimes E)$ we need to understand each $H^0(M,T\otimes {\lie g}_{x})$. The bundles ${\lie g}_{x}$ are parametrized by $x\in C$ and so in the complement of a finite set of points in $C$, $\dim H^0(M,T\otimes {\lie g}_{x})$ takes its generic value $k$, say. For any $x$ we have the section $s_x$ defined in Section \ref{poly1} and so $k\ge 1$. Since the canonical bundle has no base points, a generic canonical differential $\theta$  vanishes at  points in this complement and so for this bundle $E$ we have 
$$\dim H^0(M,T\otimes E)=\sum_{i=1}^{2g-2}\dim H^0(M,T\otimes {\lie g}_{x_i})=(2g-2)k.$$
Now $H^0(M,T\otimes T)=H^0(M,\Sym^2 T)\oplus H^0(M,\Lambda^2 T)$ and it was proved in  \cite{NJH} that  $\dim H^0(M,\Sym^2T^*)=3g-3$. Let $n$ be the dimension  of the image of $H^0(M,T\otimes T)$ in $H^1(M, T\otimes T^*)$ in the above sequence then from exactness 
$$n+(2g-2)k=1+(3g-3)+\dim H^0(M,\Lambda^2 T)$$
using again $\dim H^0(M,T\otimes T^*)=1$. But $n\ge 0$ and $\dim H^0(M,\Lambda^2T)\le 1$ so if $g>2$ we must have $k=1$ and $n=g+\dim H^0(M,\Lambda^2T)$. But $n\le \dim H^1(M,T\otimes T^*)=g$ and hence $H^0(M,\Lambda^2T)=0$.
\end{prf}

\begin{rmks}

\noindent 1. When $g=2$, $M$ is the intersection of two quadrics in the $5$-dimensional projective space $\PP(V)$. A direct calculation shows that $H^0(M,\Lambda^2T)\cong \Lambda^2V^*$.

\noindent 2. From \cite{Toda} the infinitesimal deformations of the abelian category of coherent sheaves are parametrized by the Hochschild cohomology group $H\!H^2(M)$ and the vanishing of $H^0(M,\Lambda^2T)$ and $H^2(M,{\mathcal O})$ means that this is isomorphic to $H^1(M,T)$, the deformations of the complex structure of $M$ which is well-known to be canonically isomorphic to the deformations of the curve $C$.

\noindent 3. The evaluation map $H^0(C,{\lie g}\otimes K) \rightarrow {\lie g}_x\otimes K_x$ defines as in  Section \ref{poly1} a holomorphic section $s_x$ of $\Hom(T^*, {\lie g}_x)=T\otimes {\lie g}_x$ on $M$. Our calculation above of $k=1$ shows that for generic $x$ this is the unique section. 

\end{rmks}

\section{Generators and relations}
\subsection{Generators}
Suppose now that $M$ is the moduli space of rank $2$ bundles of  fixed determinant over a curve $C$ of genus $g$. We have seen from Proposition \ref{inj} that the $(5g-5)$-dimensional space $H^1(C,K^{-2})$ injects into $H^0(M,\Lambda^3T)$. This generates maps
$$\Lambda^k  H^1(C,K^{-2})\rightarrow H^0(M,\Lambda^{3k}T)$$
and one may ask whether this is surjective, or more generally is it true that any polyvector field is generated by these trivector fields? 

Since $\dim M=3g-3$ we can consider the map from  $\Lambda^{g-1}  H^1(C,K^{-2})$ to sections of the anticanonical bundle $K_M^{-1}= \Lambda^{3g-3}T$ of $M$. The Verlinde formula gives this dimension as 
$$\dim H^0(M,K_M^{-1})=3^{g-1}2^{2g-1}\pm 2^{2g-1}+3^{g-1}$$
(where the sign corresponds to even or odd degree),
whereas
$$\dim \Lambda^{g-1}  H^1(C,K^{-2})={5g-5\choose g-1}$$
which is  smaller.

On the other hand, our polyvector fields are described via the adjoint representation and so are insensitive to the operation of tensoring a rank $n$ stable vector bundle $V$ of fixed determinant by a line bundle of order $n$. So on the moduli space $M$ of stable vector bundles  they are invariant by the action of $H^1(C,\Z_n)$. In the rank $2$ case the dimension of the space of invariant sections of $K^*$ is given in \cite{OP} as 
$$\dim H_0^0(M,K_M^{-1})=\frac{3^g\pm 1}{2}.$$
Using the inequality
$${n\choose k}\ge \left(\frac{n}{k}\right )^k$$
we have for $g>2$ 
$$\dim \Lambda^{g-1}  H^1(C,K^{-2})={5g-5\choose g-1}\ge 5^{g-1}\ge \frac{3^g\pm 1}{2}.$$
It therefore remains a possibility that the invariant trivectors do generate the whole algebra.

\subsection{Some relations}
Recall that for each point $x \in C$  we have (up to a constant) a trivector $\sigma_x$ defined by evaluation at $x$:
$$\sigma_x(\Phi_1,\Phi_2,\Phi_3)=B(\Phi_1(x),[\Phi_2(x),\Phi_3(x)]).$$
For $SL(2)$ the three-form $B(X,[Y,Z]))$ is essentially the volume form of the Killing form on the three-dimensional Lie algebra.

Now take $(g-1)$ distinct points $x_1,\dots,x_{g-1}$ on $C$ and consider evaluating a Higgs field $\Phi$, considered as a cotangent vector to $M$,  at these points. We get a homomorphism 
$$\alpha: T^*\rightarrow \bigoplus_{i=1}^{g-1}{\lie g}_{x_i}$$
of bundles of the same rank. Taking the top exterior power
$$\Lambda^{3g-3}\alpha: \Lambda^{3g-3}T^*\rightarrow \bigotimes_{i=1}^{g-1}\Lambda^3 {\lie g}_{x_i}.$$
The right hand side is just a trivial bundle so this homomorphism defines a  section of the anticanonical bundle of $M$ naturally associated to the $(g-1)$ points. In fact it is not hard to see that it is a multiple of 
$$\sigma_{x_1}\wedge\sigma_{x_2}\wedge\dots\wedge\sigma_{x_{g-1}}.$$
This vanishes when $\alpha$ has a non-zero kernel, which is the locus of bundles in $M$ for which there is a Higgs field vanishing at the $(g-1)$ points -- a determinant divisor.

If the rank $2$ vector bundle has degree zero then by the mod 2 index theorem (as for example in \cite{Bea}), if $K^{1/2}$ is an odd  theta characteristic  then 
$$\dim H^0(C,{\lie g}\otimes K^{1/2})>0.$$
So if $\Psi\in   H^0(C,{\lie g}\otimes K^{1/2})$ and a section $s$ of $K^{1/2}$ has divisor $x_1+x_2+\cdots + x_{g-1}$ then $\Phi=s\Psi$ is a Higgs field which vanishes at these points. In other words every bundle has a Higgs field vanishing at these points so
$$\sigma_{x_1}\wedge\sigma_{x_2}\wedge\dots\wedge\sigma_{x_{g-1}}=0.$$
These are relations in the algebra -- one for each  of the $2^{g-1}(2^g-1)$ odd theta characteristics. 
However we still have for $g>4$ 
$$5^{g-1}-2^{g-1}(2^g-1)>\frac{3^g+ 1}{2}$$
so there must be more.

\vskip 1cm
 Mathematical Institute, 24-29 St Giles, Oxford OX1 3LB, UK
 
 hitchin@maths.ox.ac.uk

 \end{document}